\documentclass[12pt]{article}
 \usepackage[english]{babel}
\usepackage[T2A]{fontenc}
\usepackage{amssymb}
\usepackage{amsmath}
\usepackage{amscd}
\usepackage{setspace}
\begin{document}
\newtheorem{mydef}{Definition}[section]
\newtheorem{theorem}{Theorem}[section]
\newtheorem{seqv}{Corollary}[section]
\newtheorem{lemm}{Lemma}[section]
\newtheorem{examp}{Example}[section]
\newtheorem{abst}{Abstract}[section]
\newenvironment{proof}
    {\par\noindent{\bf Proof.\enspace}}
    {\hfill$\square$}
\author{Z.R.~Gabidullina}
\title{The Minkowski Difference for Convex Polyhedra and Some its Applications
\thanks{Kazan Federal University,\, e-mail: zgabid@mail.ru,\, zulfiya.gabidullina@kpfu.ru }}
\date{}
\maketitle
 \abstract The aim of the paper is to develop a~unified algebraical
approach to representing the Minkowski difference for convex
polyhedra. Namely, there is proposed an exact analytical formulas
of the Minkowski difference for convex polyhedra with different
representations. We study  the cases when both operands under the
Minkowski difference operation simultaneously have a vertex or a
half-space representation. We also focus on the description of the
Minkowski difference for a such mixed case where the first operand
has the linear constraint structure and the second one is
expressible as the convex hull of a finite collection of some
given points. Unlike the widespread geometric approach considering
mostly two-dimensional or three-dimensional spaces, we investigate
the objects in finite-dimensional spaces of
arbitrary dimensionality.\\
{\bf keywords:} Minkowski difference, convex polyhedron, vertex
representation, half-space representation,
 polyhedra, distance, projection,  linear separability
 criterion,   variational inequality\\
 {\bf MSC classes:} 90C30,\,  65K05

\section{Introduction}\label{Gabid_sec1}

 The Minkowski difference of  sets having the different
configuration is basic to the treatment of a large class of
problems often occurring in lots of in\-te\-res\-ting applications
from a variety of areas, especially in problems of engineering
design \cite{G_Zagajac},\, \cite{G_Ilies}; in data classification
\cite{Mampaey},\, \cite{Takeda}, \cite{Mavroforakis}; image
analysis and processing \cite{G_Serra},\, \cite{G_Barki};\, motion
planning for robots \cite{G_Lozano-Perez},\,\cite{G_Cameron};\,
real-time collision detection \cite{G_Ericson};\, computer
graphics \cite{G_Ghosh1},\, \cite{G_Ghosh}, and many other
front-line fields.

The Minkowski difference operation is actually useful as an
investigative as well as a conceptual tool.
 But unfortunately, it is widely known fact that
  there exist serious difficulties related to
implementation of the Minkowski difference for individual
formulations of sets. They represent the basic impediment to
making use of the Minkowski difference operation in various
practical applications. For finite-dimen\-sional spaces of
arbitrary dimensionality, an exact analytical representation of
the Minkowski difference for the convex poly\-hed\-ra with the
different as well as identical configuration is stated here for
the first time as a whole.

To be adequate for a number of  mathematical purposes,  the
different approaches (with cross-fertilization of ideas) to such a
basic concept of analysis  as the Minkowski difference of sets are
required. The geometric viewpoint among the others has the leading
role. Motivation for the development of this geometric approach
has basically come from engineering design \cite{Nelaturi}
\cite{G_Cappelli}, computational geometry \cite{G_O'Rourke},\,
\cite{Klette},\, \cite{Karasik}, collision detection
\,\cite{G_Bergen},\, \cite{G_Bergen1}, robotics
 \cite{G_Gilbert},\, \cite{G_Gilbert1} and many
other subjects. The nature of these subjects dictates the
sufficiency of considering only two-dimensional or
three-dimensional spaces (see, for instance, \cite{G_Barki},\,
\cite{G_Ghosh},\, \cite{G_Hachenberger}).\, Unlike the geometric
approach, we study the objects in spaces of arbitrary
dimensionality and develop the algebraic approach.

Namely, we present the exact analytical representation of the
Min\-kow\-ski difference for convex polyhedra given by different
ways. More precisely, we investigate the following cases where:
\begin{itemize}
    \item both operands under the Minkowski difference operation are
determined in the similar way as the convex hull of the finite
collection of some given points,
 \item the operands have the different representation (namely, the
first operand is given by a linear constraint system, and the
second one is expressible in terms of the convex hull),
 \item both operands have the same representation as the intersection
 of the closed half-spaces.
    \end{itemize}
Let us note that thanks to the obtained results, there appeared
the possibility to treat the relevant problems such as:
\begin{itemize}
    \item problem of
linear se\-pa\-ra\-tion of the convex polyhedra in the Euclidean
space, \item variational inequalities problems, \item problem of
finding the distance between the convex polyhedra by projecting
the origin of the Euclidean space onto a convex polyhedron, \item
problem of finding the nearest points of convex polyhedra.
 \end{itemize}

\section{Definitions and Preliminaries}
This section includes
 the brief description of notations, definitions of all utilized
 in the present paper main concepts. We use standard notation that is certainly well
known to all readers. Nevertheless, let us briefly describe some
of notations. As usual, \,$\|\cdot\|$\, stands for the
Euc\-li\-de\-an norm  of the vector in $\mathbb{R}^{n}$,\,
\,$\langle \cdot, \cdot \rangle$\, denotes the standard scalar
product of two vectors, \,$conv\{\cdot\}$\, corresponds to the
convex hull of some collection of the given vectors.
  Let \,$Bd(\cdot)$,\,\,$int(\cdot)$\,
 denote the boundary and interior of some set \,$\tilde{\Phi}$,\, respectively.

 In a context of applications for the concept of the Minkowski
difference, we need further
 to recall some key definitions and theorems which are closely related
  to  the linear separability property of sets. In the theory of convex sets
  and nonlinear programming, the study of the linear separation
  problem is a~high\-ly im\-por\-tant topic  with a large
literature (see, for instance, \cite{G_Giannessi},\,
\cite{Gabid_Vas},\,  \cite{G_Eremin},\, \cite{G_Gabid2},\,
\cite{G_Gabid6}, \cite{G_Gabid},\, \cite{G_Gabid7}, etc.).
  \begin{mydef}\label{G_d1} (separating hyperplane) (see, for instance, \cite{Gabid_Vas}, p.
          198)\quad
The hyperplane \smallskip \\
\centerline{$\pi(c,\gamma): = \bigr \{x \in~
\mathbb{R}^{n}:~\langle c,x \rangle =\gamma \bigl \}$}
\\ with normal vector $\,c \neq \bf 0 \,$\, separates
the sets $A$ and $B$ from the Euclidean space $\mathbb{R}^{n}$,
 iff $\,\langle c,a\rangle \geq \gamma\,$ for all
$\,a\in A\,$ and $\,\langle c,b \rangle \leq \gamma \,$ for all
$\,b\in B$,\, i.e., iff it holds that:\\
\centerline {$\sup \limits_{b\in B}\,\, \langle c,b \rangle \leq
\gamma\leq \inf \limits _{a \in A}\,\, \langle c,a \rangle.$}
        \end{mydef}
    \begin{mydef}\label{G_d2}\,(strong separability) (see \cite{Gabid_Vas}, p. 198)\,
Two sets \,$A$ and $B$\, are said to be strongly separable, iff
there exists some vector $\,c \in \mathbb{R}^{n}$\, such that:\\
\centerline{$\sup \limits_{b\in B}\,\, \langle
c,b \rangle < \inf \limits _{a\in A}\,\, \langle c,a \rangle.$}\\
 If it holds that \,$\langle c,b \rangle < \langle c,a \rangle$\,
for all \,$a \in A,$\, \,$b \in B$,\, then the sets \,$A$\, and
\,$B$\, are said to be strictly separable.
    \end{mydef}
 The next theorem gives the rigorous justification of
the fact that the problem of strong separation of the two
arbitrary sets \,$A$,\, $B \subset \mathbb{R}^{n}$\, can be
reformulated as the problem of strong separating the origin of
\,$\mathbb{R}^{n}$\, from the Minkowski difference \,$A - B = \{a
- b, a \in A,\, b \in B\}$,\, and vice versa.
\begin{theorem} \label{G_2.7} (strong separation) (\cite{G_Gabid6},\, p. 150)
For the sets \,$A$ and \,$B$\, to be strongly separable, it is
necessary and sufficient that the origin of \,$\mathbb{R}^{n}$ be
strongly separable from the set \,$A-B.$
         \end{theorem}
         The previous theorem immediately implies that two sets \,$A$\, and
\,$B$\, are not strongly separable if and only if the set
\,$A-B$\, is not strongly separable from the origin of
\,$\mathbb{R}^{n}$.\,

         The analogous results are certainly well known for the problems
          on non-strong linear separation of the considered sets.
         \begin{theorem} \label{G_2.77} (linear separation)(\cite{G_Gabid6},\, p. 151)
For two sets \,$A$\, and \,$B$\, to be linearly separable, it is
necessary and sufficient that the origin of \,$\mathbb{R}^{n}$ be
linearly separable from the set \,$A-B.$\,
         \end{theorem}
         The preceding theorem obviously implies that
to say two sets \,$A$\, and \,$B$\, are linearly inseparable is to
say the origin of \,$\mathbb{R}^{n}$ is linearly inseparable from
the set \,$A-B.$\,

 We note that
 if the sets \,$A$ and \,$B$\, are convex, then the Minkowski
 difference \,$A-B$\, is convex as well (see, for example, \cite{Gabid_Vas}, p. 162).
 It is not hard to prove that if \,$A$ and \,$B$\,
 are simultaneously bounded sets, then the set \,$A-B$\, is also bounded.
 Lastly, under the condition that at least one of the sets \,$A$ and \,$B$\, is bounded,
 closedness of both sets under the Minkowski difference operation implies
  closedness of \,$A-B$.\, The proof of this assertion may be found, for instance,
   in \cite{Gabid_Vas}, p. 201.

 Another question can now be addressed. What conditions on some
sets \,$A$,\, $B$\, ensure that they are linearly (strongly or
not) separable from each other?  As already noted above, we reduce
the problem of linear separation of \,$A$\, and $B$\, to the
program of separation the origin of $\mathbb{R}^{n}$\, from \,$A -
B$.\, Questions about the separability tests or criteria for sets
can actually be answered from the platform of the Minkowski
difference operation. It is widely known that the Minkowski
difference is often interpreted as the translational configuration
space obstacle (see, for instance, \cite{G_Cameron}).   One says
that \,$A - B$\, represents the set of translations of \,$B$\,
that brings it into interference with \,$A$.\, The additional nice
property of the Minkowski difference consists in that for any
objects \,$A$\, and $B$\, it holds
$$dist(A,\,B) = dist({\bf 0},\,A - B),$$
$$\text{where} \enspace dist(A,\,B) = \inf \{\|a - b\|,\, a \in A,\, b \in B\}$$
denotes the distance between $A$ and $B$. Clearly, \,$dist({\bf
0},\,A - B) = \|{\bf P}_{A - B}({\bf 0})\|,$\, where \,${\bf P}_{A
- B}({\bf 0})$\, denotes the projection of the origin of
\,$\mathbb{R}^{n}$\, onto \,$A - B$.\, As it is known, two convex
objects collide if and only if their Minkowski difference contains
the origin.
 For
the origin of $\mathbb{R}^{n}$\, and \,$A - B$,\, the next results
of this subsection give a sort of linear separation principle.

 First we recall briefly
some other relevant notations and definitions about the
 cones of generalized strong and strict support
vectors, etc. In \cite{G_Gabid6},\, there have been defined the
following sets:
 $$W_{\tilde{\Phi}} := \bigl \{w \in \mathbb{R}^{n}:
 \enspace \inf \limits_{x \in \tilde{\Phi}}\, \langle w, x \rangle \geq 0 \bigr\},
 \,V_{\tilde{\Phi}} := \bigl \{v \in \mathbb{R}^{n}:
  \enspace \inf \limits_{x \in \tilde{\Phi}}\, \langle v, x \rangle > 0 \bigr\},$$
  $$Q_{\tilde{\Phi}} := \bigl \{q \in \mathbb{R}^{n}:
  \enspace \langle q, x \rangle > 0, \enspace \forall x \in \tilde{\Phi} \bigr\},$$
$$\Omega_{\tilde{\Phi}} := \{y \in \mathbb{R}^{n}:
\,\langle y,x \rangle \geq \|y\|^{2}, \enspace \forall x \in
\tilde{\Phi}\},$$
  $$E(\Omega_{\tilde{\Phi}}):= \{ x \in \mathbb{R}^{n}: x = \lambda y, \lambda \geq 0,\,
  y \in \Omega_{\tilde{\Phi}}\},$$
 where \,$\tilde{\Phi}$\, is a nonempty subset
of \,$\mathbb{R}^{n}$.\, The set \,$W_{\tilde{\Phi}}\backslash
\{{\bf 0}\}$\, is called a cone of generalized support vectors
(or, briefly, GSVs) of the set \,$\tilde{\Phi}$.\, The notations
\,$V_{\tilde{\Phi}}$\, and \,$Q_{\tilde{\Phi}}$\, are used for the
cones of generalized strong and strict support vectors of the set
\,$\tilde{\Phi}$,\, respectively. The main properties of the
mentioned cones and their relationship with well-known other ones
have been investigated in \cite{G_Gabid6}.\, For details, the
interested reader is recommended to refer to \cite{G_Gabid6}.
However, recalling some main properties of GSVs seems as vital for
understanding  the results of the present paper:
$$W_{\tilde{\Phi}} = W_{c\,l(\tilde{\Phi})} = W_{co(c\,l(\tilde{\Phi}))},
\enspace V_{\tilde{\Phi}} = V_{co(\tilde{\Phi})} =
V_{c\,l(co(\tilde{\Phi}))},$$
$$V_{\tilde{\Phi}}\subseteq
Q_{\tilde{\Phi}}\subseteq W_{\tilde{\Phi}}, V_{\tilde{\Phi}} =
E(\Omega_{\tilde{\Phi}})\backslash\{{\bf 0}\},$$
$$V_{\tilde{\Phi}} = V_{c\,l(\tilde{\Phi})} = V_{co(c\,l(\tilde{\Phi}))},
 \enspace W_{\tilde{\Phi}} = W_{co(\tilde{\Phi})} =
W_{c\,l(co(\tilde{\Phi}))}, \enspace Q_{\tilde{\Phi}} =
Q_{co(\tilde{\Phi})}\,.$$
 The conditions which are necessary and
sufficient for the emptiness of the cones of GSVs
\,$V_{\tilde{\Phi}}$,\, \,$Q_{\tilde{\Phi}}$,\, and
\,$W_{\tilde{\Phi}}\backslash \{{\bf 0}\}$\, were established in
\cite{G_Gabid7}.

Let, \,$t_{\tilde{\Phi}}(c) := \inf \limits_{x \in
\tilde{\Phi}}\,\,\langle c, x \rangle$,\, \,$\tilde{\Phi} \subset
\mathbb{R}^{n}$.\, For some kinds of convex polyhedra which have
the compactness property,  due to the presence of convexity and
continuity of the linear function \,$\langle c,x \rangle$,\,  we
can state that \,$\langle c,x \rangle$\, attains its minimum on
\,$\tilde{\Phi}$.\, So, for these cases, we can rewrite
\,$t_{\tilde{\Phi}}(c) := \min \limits_{x \in
\tilde{\Phi}}\,\,\langle c, x \rangle$.\, Further, we note that
the following problem
                  \begin{equation}\label{G_maxtfc}
                  \max \limits_{\|c \| = 1}
                  \,\,t_{\tilde{\Phi}}(c)
                   \end{equation}
is solvable. For the proof of the fact, the interested reader is
directed to \cite{G_Gabid6}\, (see p. 149). Let the vector
\,$c^{*}\in \mathbb{R}^{n}$\, denote an~optimizer of the problem
(\ref{G_maxtfc}). The following three theorems describe a linear
separability criterion for the pair of objects such as the origin
and some nonempty set of \,$\mathbb{R}^{n}$.\, Based on the
optimal value of the objective function of~(\ref{G_maxtfc}),\, the
above-mentioned criterion allows us to recognize these objects as
strongly separable, non-strongly linearly separable, or
inseparable.
                     \begin{theorem} \label{G_2.9} (strong separability criterion)
                     (\cite{G_Gabid6},\,\,p.149)
\,For the origin of \,$\mathbb{R}^{n}$ to be strongly separable
from the nonempty set \,$\tilde{\Phi} \subset \mathbb{R}^{n}$\, it
is necessary and sufficient to have \,\,$t_{\tilde{\Phi}}(c^{*}) >
0$.\,
                     \end{theorem}
                     \begin{theorem} \label{G_2.8}(non-strong linear separability criterion)
                     (\cite{G_Gabid6},\,p.150)
    For the origin of \,$\mathbb{R}^{n}$\, to be non-strongly linearly
    separable from the nonempty set \,$\tilde{\Phi}~\subset~\mathbb{R}^{n}$\, it is
    necessary and sufficient to have \,\,$t_{\tilde{\Phi}}(c^{*}) = 0$.\,
                        \end{theorem}
                        \begin{theorem} \label{G_2.10} (linear inseparability criterion)
                        (\cite{G_Gabid6},\,p.150)
\,For the origin of \,$\mathbb{R}^{n}$ to be linearly inseparable
from the nonempty set \,$\tilde{\Phi}\subset~\mathbb{R}^{n}$\, it
is necessary and sufficient to have \,\,$t_{\tilde{\Phi}}(c^{*}) <
0$.\,
                        \end{theorem}
                        The following theorems allows to detect which one of~the
cones of generalized support vectors (\,$V_{\tilde{\Phi}}$,
\,$W_{\tilde{\Phi}}\backslash \{{\bf 0}\}$,\, and
\,$Q_{\tilde{\Phi}}$) is empty, and which ones are not.
                         \begin{theorem}\label{G_concl2} (emptiness of the cone
                          of generalized strong support vectors)
                           If  \,$\tilde{\Phi}$\, is a nonempty
                convex and closed subset of  \,$\mathbb{R}^{n}$,\,  then $$V_{\tilde{\Phi}}
= \emptyset \,\,
 \,\,\Leftrightarrow \,\,{\bf 0}\in \tilde{\Phi}.$$
        \end{theorem}
        The previous theorem represents the particular case of Theorem 3.3
        from \cite{G_Gabid7} (for the proof, see p. 703).\, For the convex and closed
         set \,$\tilde{\Phi}$,\,
         due to having \,$V_{\tilde{\Phi}} = E(\Omega_{\tilde{\Phi}})\backslash\{{\bf 0}\}$,\,
         from Theo\-rem~\ref{G_concl2},\,
        there obviously holds the following implication: \,$\Omega_{\tilde{\Phi}}
= \{{\bf 0}\} \,\,  \,\,\Leftrightarrow \,\,{\bf 0}\in
\tilde{\Phi}.$\,
         \begin{theorem}\label{G_concl1} (emptiness of the cone of generalized
         support vectors)\\ If \,$\tilde{\Phi} \subset \mathbb{R}^{n}$\, is a nonempty
                convex set,
then \,$$W_{\tilde{\Phi}} = \{{\bf 0}\}\,\, \,\Leftrightarrow\,\,
\,\,{\bf 0} \in int(\tilde{\Phi}).$$
       \end{theorem}
       The preceding theorem is the particular case of Theorem 3.2
        from \cite{G_Gabid7} (for the proof, see p. 703).\,
              \begin{theorem}\label{G_Utv_19} (simultaneous degeneracy of the cone of generalized
              strict support vectors $\&$ non-degeneracy of the cone of GSVs)
                If \,$\tilde{\Phi}\subset \mathbb{R}^{n}$\,  is a nonempty convex and closed set,
                then:\\
                \centerline{$\,W_{\tilde{\Phi}} \neq \{{\bf 0}\}$\, $\&$ \,$Q_{\tilde{\Phi}} =
                 \emptyset$\, $\Leftrightarrow$
                 ${\bf 0} \in Bd(\tilde{\Phi})$.}
       \end{theorem}
       According to \cite{G_Gabid7}, the assertion of the preceding theorem follows from
       Lem\-mas~3.13--3.14 (see p. 702).\, Let us note that for
       the set \,$\tilde{\Phi}$\, having the compactness property, it is
       fulfilled \,$V_{\tilde{\Phi}} = Q_{\tilde{\Phi}}$.\,
       According to Theorems \ref{G_2.9},\, \ref{G_concl2},\, we
       obviously have $${\bf 0}\notin \tilde{\Phi} \,\,\Leftrightarrow\,\, t_{\tilde{\Phi}}(c^{*}) >
0.$$ From Theorems \ref{G_2.10},\, \ref{G_concl1},\, it
immediately follows that $${\bf 0} \in int(\tilde{\Phi})
\,\,\Leftrightarrow\,\, t_{\tilde{\Phi}}(c^{*}) < 0.$$ Due to
Theorems \ref{G_2.8},\, \ref{G_Utv_19},\, there holds the
following implication: $${\bf 0} \in Bd(\tilde{\Phi})
\,\,\Leftrightarrow\,\, t_{\tilde{\Phi}}(c^{*}) = 0.$$ Thus, for
the origin of \,$\mathbb{R}^{n}$,\, the linear separability
criterion  provides a certificate of being an~exterior, interior,
or boundary point of \,$\tilde{\Phi}$.\,
\section{Binary Operation of Minkowski Difference}\label{Gabid_sec2}
\subsection{Both Operands with a Vertex Representation}\label{Gabid_sbsec21}\,
For many applications nowadays, the sets expressible as the convex
hull of finitely many points from the Euclidean space
\,$\mathbb{R}^{n}$\, are especially important.  They are really
ubiquitous structures having a fundamental role not only in
variational analysis, computational geometry and optimization, but
in data classification, image analysis and processing, motion
planning for robots, collision detection, and many other
front-line areas. This subsection is devoted to representing the
Minkowski difference for both convex polyhedra having the
above-mentioned configuration.

 Further, let
us be given two polyhedra of \,$\mathbb{R}^{n}:$\,
\begin{multline}\nonumber
L:= conv \{z_{i}\}_{i \in
  I} \,M:= conv \{p_{j}\}_{j \in J}, \smallskip \\ \text{where}
  \,I = \{1,2,\cdots, l \},\,
  \,J = \{1,2,\cdots, m \},\, \text{i.e.}
\end{multline}
$$L = \{ z\in \mathbb{R}^{n}: z=\sum \limits_{i \in I} \alpha_{i}z_{i}, \,
 \, \sum \limits_{i \in I}\alpha_{i}=1, \,\alpha_{i}\geq 0, \,i \in I \},$$
$$M = \{p \in \mathbb{R}^{n}: p=\sum \limits_{i \in J} \beta_{j}p_{j},
 \, \sum \limits_{j \in J}\beta_{j}=1, \, \beta_{j}\geq 0, j \in J \}.$$

Obviously, \,$L$ and \,$M$\, are nonempty, convex, and compact
sets. Due to set algebra, the Minkowski difference of these
polyhedra is defined as a set of pairwise differences of points
from \,$L$ and \,$M$.\, Namely, as follows \smallskip \\
\centerline{ \,$L-M = \{z-p, \,z \in L,p \in M \}$.\,}\smallskip

 The following theorem  rigorously justifies that the point set \,$L-M$\,
 coincides with the convex hull
 of the vectors  \,$z_{i}-p_{j},$ \enspace $i \in I,$ $j \in J$.
\begin{theorem}\label{G_t2}. (Minkowski difference for both polyhedra given as convex hull)  \,(\cite{G_Gabid2},\, p. 552)\,
\,$L-M = conv \{z_{i}-p_{j}\}_{i \in I, j \in J}.$\,
\end{theorem}
\proof Part I.\, At first, we establish the inclusion\\
\centerline{\,$conv \{z_{i}-p_{j}\}_{i \in I, j \in J} \subseteq
L-M.$}\\ By the definition of the convex hull, for any \,$l \in
conv \{z_{i}-p_{j}\}_{i \in I, j \in J}$,\, there can be found the
real numbers \,$\gamma_{ij}\geq 0,$ \,$\sum \limits_ {i \in I}\sum
\limits_ {j \in J} \gamma_{ij}=1$\, such that: $$l=\sum \limits_
{i \in I}\sum \limits_{j \in J} \gamma_{ij}(z_{i}-p_{j})= \sum
\limits_{i\in I} (\sum \limits_{j \in J} \gamma_{ij})z_{i}- \sum
\limits_{j\in J}(\sum \limits_{i \in I} \gamma_{ij})p_{j}.$$
Assuming from now that \,$\alpha_{i}=\sum \limits_{j \in
J}\gamma_{ij}, \beta_{j}=\sum \limits_{i \in I}\gamma_{ij},$\, one
can easily see that the coefficients \,$\alpha_{i}$ and
$\beta_{j}$\, satisfy the conditions \,$\alpha_{i}\geq 0,\,
\forall i \in I,$\, $\sum \limits_{i \in I}\alpha_{i}=1,$\,
\,$\beta_{j}\geq 0,$ $\forall j \in J,$\, \,$\sum \limits_{j \in
J}\beta{j}=1.$\, Consequently, any vector \,$l \in conv
\{z_{i}-p_{j}\}_{i \in I, j \in J}$\, satisfies the following
equation:\\

 \centerline{$l=\sum \limits_{i \in I}
\alpha_{i}z_{i}-\sum \limits_{j \in J}\beta_{j}p_{j}=z-p,\,z \in
L,\, p \in M,$ \text{i.e.} $l \in L-M.$}
Part II.\, Now let us justify the following backward inclusion:\smallskip\\
\centerline{$L-M \subseteq conv \{z_{i}-p_{j}\}_{i \in I, j \in
J}.$}\smallskip \\ By the definition of the set \,$L-M$,\, taking
 some vector \,$l \in L-M,$\,
 we then have that \,$l=z-p, \, z
\in L, p \in M.$ According to the construction of the sets $L$ and
$M$, there will be found \,$\alpha_{i}\geq 0$ \,$\forall i \in I,$
$\beta_{j}\geq 0$ \enspace $\forall j \in J,$ \,$\sum \limits_{i
\in I}\alpha_{i}=1,$ \,$\sum \limits_{j \in J}\beta_{j}=1$\, such
that \,$z=\sum \limits_{i \in I} \alpha_{i}z_{i},\, p=\sum
\limits_{j \in J}\beta_{j}p_{j}.$ Then
\begin{multline}\nonumber
  l=\sum \limits_{i \in I} \alpha_{i}z_{i}- \sum \limits_{j \in
J}\beta_{j}p_{j}=\sum \limits_{j \in J}\beta_{j}(\sum \limits_{i
\in I} \alpha_{i}z_{i})- \sum \limits_{i \in I} \alpha_{i}(\sum
\limits_{j \in J}\beta_{j}p_{j}) = \\ = \sum \limits_{i \in I}\sum
\limits_{j \in J} \alpha_{i}\beta_{j}(z_{i}-p_{j})= \sum
\limits_{i \in I}\sum \limits_{j \in J} \gamma_{ij}(z_{i}-p_{j}),
\end{multline}
where \,$\gamma_{ij}\geq 0,$ \,$\sum \limits_{i \in I}\sum
\limits_{j \in J} \gamma_{ij}=1.$ Thus, there is true the
following inclusion:  \,$L-M \subseteq conv \{z_{i}-p_{j}\}_{i \in
I, j \in J}.$\, \smallskip The latter should be compared to the
earlier proved inclusion \,$conv \{z_{i}-p_{j}\}_{i \in I, j \in
J} \subseteq L-M$.\, This comparison allows us to complete the
proof.  $\square$ \\ Since due to Theorem~\ref{G_t2},\, \,$L-M$\,
has a representation as a convex hull of fi\-ni\-te number
\,$(l\cdot m)$\, of points \,$z_{i}-p_{j}$\, for all \,$i \in I,$
 $j \in J$,\, it is characterized as a nonempty, compact, and convex
 point set. Consequently, for the case, the operation of the Minkowski
 difference thereby preserves the compactness and
 convexity.

 For practical applications, there has the extremely importance a question
 consisting of how to decrease the number of points forming \,$L-M$.\,
 For the low enough dimension of \,$\mathbb{R}^{n}$,\, we have a possibility
 to make this number of points as small as possible by means of using some software package.
 In the case of the two-dimensional or three-dimensional space, the function Convhull,
 for instance, in MATLAB  not only computes and  returns the convex hull of the
given collection of points, but provides the option of removing
vertices that do not contribute to the area or volume of the
convex hull. Moreover, this package allows to visualize the output
of Convhull with the help of the function Plot in 2-D. The
function Trisurf or Trimesh provides the possibility of plotting
the output of Convhull in 3-D. Let us note that in four or more
dimensions, there can efficiently be used, for instance, the
proposed in \cite{G_Barber} Quickhull algorithm for computing the
convex hulls.  This method is realized in MATLAB by means of the
function Convhulln. This function returns the indices of input
points that form the faces of the convex hull. Consequently, to
compute, for instance, the distance between \,$L$ and \,$M$\, or
to linearly separate these sets, one first should select those of
the \,$(l\cdot m)$\, points of the type \,$z_{i}-p_{j}$,\, \,$i
\in I,$ $j \in J$\, that are really formed the faces of \,$L-M$.\,
Luckily, the inner points of the Minkowski difference \,$L-M$\,
may be ignored. As a result of taking into account of the only
points that are vertices of the convex hull, the collection of
points \,$z_{i}-p_{j}$,\, \,$i \in I,$ $j \in J$\, may
considerably be reduced. Therefore, it is more easier to deal with
such a reduced family of input points. The expected time
complexity of the Quickhull algorithm depends on the different
parameters such as dimension of the space, the number of input and
processed points, the maximum number of facets (for details, see
\cite{G_Barber}, \cite{G_David}).
\subsection{Binary Mixture of Sets with Different Constraint Structure}\label{Gabid_subsec2.2}\,
In this subsection, we focus on a precise representation of the
Minkowski difference for a common case where the first point set
from the sets pair under the operation has the general constraint
(i.e. not necessarily linear constraint) structure, and the second
operand is identified by the abstract constraint. Such sets
identification is exhibited as too crucial for applications to be
considered below. Furthermore, the abstract constraint
specification is a~real\-ly useable form since it does not
restrict the variations on how the corresponding point set might
be defined. The presence of the abstract constraint characterizes
our approach as very flexible \cite{G_Gabid7},\,\cite{G_Gabid10},
since constraints might not even be present. The purpose of
considering in this subsection such a case of the more general
settings is twofold - to recall some fundamental results from the
previous research and to explain how they can be applied to
a~to\-pic of our interest.

  In a wide range of practical
applications, the set of feasible solutions \,$\Phi$\, is
representable by a system of inequality constraints in the general
form:
 \begin{equation} \label{G_sys_fhi}
 \Phi = \{x \in X : f_{k}(x) \leq b_{k}, \, k \in
K\}, \enspace K = \{1,2,\cdots , r\},
         \end{equation}
  where \,$f_{k}(x), \, k
\in K$\, are arbitrary real-scaled quasi-convex functions
 defined on a convex set
 \,$X \subseteq \mathbb{R}^{n}$.\,
We recall that a function \,$f(x)$\, is said to be
a~qua\-si-con\-vex on a~con\-vex set \,$X$\, if and only if
\,$[S^{d},f]^{Lo}_{X}$\, is a convex set for all \,$d \in
\mathbb{R}^{1}$,\, where
$$[S^{d},f]^{Lo}_{X} := \{ x \in X: f(x) \leq d \}.$$
Therefore, as an intersection of the convex sets
\,$[S^{b_{i}},f_{k}]^{Lo}_{X},\, k \in K$,\,  the set \,$\Phi$\,
is convex, too.
  \begin{theorem}\label{G___lower_semicont} (closedness of the lower level set)
        Let \,$X$\, be a closed set of \,$\mathbb{R}^{n}$,\, then
        \,$f(x)$\, is a lower semicontinuous  function over \,$X$\, if and
        only if  the lower level set \,$[S^{d},f]^{Lo}_{X}$\, is closed for
         all \,$d \in \mathbb{R}^{1}$.\,
            \end{theorem}
The interested reader can find the proof of the theorem, for
instance, in \cite{Gabid_Vas}\, (see p. 81).

 The previous theorem implies that
if \,$f_{k}(x), \, k \in K$\, are lower semicontinuous  functions
 over a closed set \,$X$,\,  then the lower level sets \,$[S^{b_{k}},f_{k}]^{Lo}_{X}$\,
  are closed
for all \,$b_{k},\, k \in K$.\, Consequently, being an
intersection of closed sets \,$[S^{b_{k}},f_{k}]^{Lo}_{X},$\, $k
\in K$,\, the set \,$\Phi$\, is closed as well.

Next, we present the analytical description of the Minkowski
difference of~two sets \,$\Phi$\, and \,$\Psi$,\, when \,$\Phi$\,
is given by (\ref{G_sys_fhi}),\, and \,$\Psi$\, is an arbitrarily
defined set.
              \begin{theorem}\label{G_Minkowski_1} (Minkowski difference when the two
              sets under the operation are given by a constraints system
              and an abstract constraint, respectively)
                 (\cite{G_Gabid7},\,\, p. 716)\,
Let be given an arbitrary  nonempty set \,$\Psi \subseteq
\mathbb{R}^{n}$,\, the set \,$\Phi \neq \emptyset$\, be defined by
(\ref{G_sys_fhi}),\,
 \,$X = \mathbb{R}^{n}$,\, then \,$\Phi
- \Psi = \Phi_{1}$,\, where
$$\Phi_{1} = \{x \in \mathbb{R}^{n}: f_{k}(x + y) \leq b_{k}, \, k \in
K,\, y \in \Psi\},$$ $\enspace K = \{1,2,\cdots , r\},$\, $\Phi -
\Psi = \{ z \in \mathbb{R}^{n}: z = x-y,\, x \in \Phi,\, y \in
\Psi\}.$
                  \end{theorem}
\proof \,Part I.\, First we select arbitrarily some fixed point
\,$\bar{x}$\, from \,$\Phi$.\, We further consider \,$\tilde{x}(y)
= \bar{x} - y$\, for \,$\forall y \in \Psi$.\, It is clear that
\,$f_{k}(\tilde{x}(y) + y) = f_{k}(\bar{x}) \leq b_{k},$\,
\,$\forall k \in K$,\, \,$\forall y \in \Psi$,\, i.e.
\,$\tilde{x}(y) \in \Phi_{1}$.\, Therefore, we have \smallskip\\
\centerline{\,$\bar{x} \in \Phi$\, $\Rightarrow$ \,$\tilde{x}(y) =
\bar{x} - y \in \Phi_{1}, \, \forall y \in \Psi$.\,}\smallskip \\
Through the fact that \,$\bar{x} \in \Phi$\, was chosen
arbitrarily,\, there was proved the following inclusion: \,$\Phi -
\Psi \subseteq \Phi_{1}$.\,

\,Part II.\, Conversely, we take now an arbitrary fixed point
\,$\bar{t} \in \Phi_{1}$\, and, for all \,$y \in \Psi$,\, check
whether the points \,$\widehat{t}(y) = \bar{t} + y$\, belong to
\,$\Phi$.\, Then there can be observed that
\,$f_{k}(\widehat{t}(y)) = f_{k}(\bar{t} + y) \leq b_{k},$\,
$\forall k \in K$, \,$\forall y \in \Psi$,\, i.e.
\,$\widehat{t}(y) \in \Phi$.\, In other words, it holds \,$\bar{t}
\in \Phi -y$,\, \,$\forall y \in \Psi$.\, Thanks to the arbitrary
choice of \,$\bar{t} \in \Phi_{1}$,\, we get that \,$\Phi_{1}
\subseteq \Phi - \Psi$.\,

We can now finalize our proof of the theorem. Since, taking into
account the forward and backward inclusions, we have the claimed
equality: \,$\Phi_{1} = \Phi - \Psi$.\,$\square$

An abstract constraint \,$y \in \Psi$\, is very convenient as well
for representing the Minkowski difference in the case of a more
complicated nature of the second point set from the sets pair
under the operation.

The previous theorem was presented for the first time in
\cite{G_Gabid7}, but it is brought out here along with the prime
role it plays in our further research.

\smallskip
In actual practice, strict inequalities are rarely seen in
constraints, however, if the nonempty  set \,$\Phi$\, happens to
be described by the strict inequality constraints,\, then
\,$\Phi_{1}$\, should be also expressed by the system of strict
in\-equali\-ti\-es.

For the proof of the previous theorem, no matter how the set
\,$\Psi$\, was expressed analytically or alternatively in other
ways. For instance, \,$\Psi$\, may be specified in a similar way
as the set \,$\Phi$:
$$\Psi := \{x \in \mathbb{R}^{n}: g_{s}(x) \leq d_{s}, \, s \in
S\},\, S = \{1, 2, \cdots , t\}.$$

 In particular, \,$\Psi$\, can consist of a single point as in the
conditions of the next lemma.
                  \begin{lemm}\label{G_Minkowski} (Minkowski difference for a set
                  with constraint structure
              and a~sing\-le\-ton)
Let be given an arbitrary vector \,$p \in \mathbb{R}^{n}$,\,
\,$\Phi$\, be defined by (\ref{G_sys_fhi}),\, \,$\Phi \neq
\emptyset$,\,
$$\Phi_{1} = \{x \in \mathbb{R}^{n}: f_{k}(x + p) \leq b_{k}, \, k \in
K\}, \enspace K = \{1,2,\cdots , r\},$$ then \,$\Phi - p =
\Phi_{1}$,\, where  $\Phi - p = \{ z \in \mathbb{R}^{n}: z =
x-p,\, x \in \Phi\}.$
                  \end{lemm}
The result of the previous lemma evidently follows from Theorem
\ref{G_Minkowski_1}\, under the assumption that \,$\Psi$\, is
a~sing\-le\-ton, i.e. \,$\Psi = \{p\}.$\,
Lemma~\ref{G_Minkowski}\, has in turn the following quite obvious
corollaries.
              \begin{seqv}\label{G_Minkowski_Positive_ort} (Minkowski difference for
              the nonnegative orthant and a singleton)
Let be given an arbitrary vector \,$p =(p^{1},\, \cdots, \,
p^{n})$,\, $\Phi$\, be described as the nonnegative orthant
$$\Phi = \mathbb{R}^{n}_{+} =  \{x =(x^{1},\, \cdots, \,
x^{n}): x^{j} \geq 0, \,j = \overline{1,\, n} \},$$ then \,$\Phi -
p = \{x =(x^{1},\, \cdots, \, x^{n}): x^{j} \geq -p^{j}, \, j =
\overline{1,\, n}\}.$\,
             \end{seqv}
            \begin{seqv}\label{G_Minkowski_Polyhedra} (Minkowski difference for
             a set with li\-ne\-ar constraint structure
              and a singleton)
Let be given an arbitrary vector \,$p \in \mathbb{R}^{n}$,\,
$\,\Phi \neq \emptyset,$
\begin{multline}
          \Phi = \{x \in \mathbb{R}^{n}: \langle a_{k},\, x \rangle \leq b_{k}, \,a_{k}
 \in \mathbb{R}^{n},\,
 b_{k} \in \mathbb{R}^{1},\, k \in
K\},\\ K = \{1,2,\cdots , r\},
\end{multline}
          then
\,$\Phi - p = \{x \in \mathbb{R}^{n}: \langle a_{k},\, x \rangle
\leq \tilde{b}_{k},\,
  \tilde{b}_{k} = b_{k} - \langle a_{k},\, p \rangle,\, k \in
K\}.$\,
             \end{seqv}
 \begin{seqv}\label{G_Minkowski_Cube} (Minkowski difference for a set with box constraints
  structure and a singleton)
Let be given an arbitrary vector \,$p \in \mathbb{R}^{n}$,\,
$\Phi$\, be specified by the box constraints on \,$x$\, of the
form
$$\Phi = \{x \in \mathbb{R}^{n}: l \leq x \leq  u, \enspace l,\, u \in
\mathbb{R}^{n}\},  \enspace \Phi \neq \emptyset,$$  then \,$\Phi -
p = \{x \in \mathbb{R}^{n}: l - p \leq x \leq u - p\}.$\,
             \end{seqv}
 \begin{seqv}\label{G_Minkowski_Ball} (Minkowski difference for
                   a closed ball and a singleton)
Let be given an arbitrary vector \,$p \in \mathbb{R}^{n}$,\,
$\Phi$\, be given as a~closed ball of radius \,$q$\, around some
point of \,$o$,\, i.e.
$$\Phi = \{x \in \mathbb{R}^{n}: \| x - o \|^{2}  \leq q^{2},
\,o \in \mathbb{R}^{n},\, q \in \mathbb{R}^{1}_{+}\},$$ then
 \,$\Phi - p = \{x \in \mathbb{R}^{n}: \| x - \bar{o}\|^{2}
  \leq q^{2}, \,\bar{o} = o - p \}.$\,
             \end{seqv}
             \subsection{Operands: Convex Polyhedra with Different Representation}\label{Gabid_subsec2.3}\,
Based on the foregoing results, we introduce in this subsection
the representation of the Minkowski difference for some binary
mixture of convex polyhedra having differently defined shapes. In
more detail, we focus now on the mixed case where the first convex
polyhedron has the linear constraint structure and the second one
is expressible as the convex hull of a finite collection of points
from \,$\mathbb{R}^{n}$.\, In this case, for the space of
arbitrary dimensionality, there is successfully reached an exact
representation of the Minkowski difference without the necessity
 of any transition to higher dimensions. Furthermore, it will be shown below that
 a number of linear constraints describing the Minkowski difference of sets
 is exactly the same as it is for the first operand.

  A set \,$\Phi \subset \mathbb{R}^{n}$\, is said to be a polyhedral set if it
can be specified  as the intersection of a finite family of closed
half-spaces, or equivalently, can be expressed by finitely many
linear constraints of form:
\begin{equation}\label{Gabid_Phi_lin_constr}
   \Phi = \{ x \in \mathbb{R}^{n}: Ax \leq b\},
\end{equation}
where \,$A$\, is the given nonvacuous matrix in \,$\mathbb{R}^{r
\times n}$\, with components \,$a_{ik}$,\, \,$b$\, is the vector
in \,$\mathbb{R}^{r}$\, with components $b_{i}$.\, It is well
known that a polyhedral set is a convex polyhedron. Let \,$M$\, be
the set of all convex combinations of some vectors \,$p_{j}, j \in
J$,\, i.e. \,$M$\, be specified in the similar manner as it was
defined in Sub\-sec\-tion~2.1.\, To construct a refined
representation of the Minkowski difference for $\Phi$\, and
\,$M$,\, we utilize first Theorem \ref{G_Minkowski_1} as follows:
$$\Phi - M = \{ x \in \mathbb{R}^{n}: A(x + y) \leq b,\, y \in M\}.$$
For the further analysis, we will need to construct the following
supplementary set:
              \begin{equation}\label{Gabid_Psi1}
   \Psi_{1} = \{ x \in \mathbb{R}^{n}: Ax \leq b - Ap_{j},\, p_{j}
 \in \mathbb{R}^{n},\, j \in J\}.
              \end{equation}
           \begin{theorem}\label{Gabid_Phi-M} (Representation of Minkowski difference
           for the mixed case before refinement)
   Let be given an arbitrary collection of vectors \,$p_{j} \in
   \mathbb{R}^{n}$,\, $j \in J,$
 \,$\Phi$\, be defined by (\ref{Gabid_Phi_lin_constr}),\, \,$\Phi
\neq \emptyset$,\, \,$M = conv\{p_{j}\}_{j \in J}$,\, then
$$\Phi - M = \Psi_{1}.$$
           \end{theorem}
\proof \,Part I.\, Without loss of generality, we fix some
arbitrary point \,$\bar{x}$\, from  \,$\Psi_{1}.$\, By our
construction of \,$\Psi_{1},$\, it then holds \,$A\bar{x} \leq b -
Ap_{j}.$\, Multiplying this system through by the nonnegative
coefficients \,$\beta_{j},$\, \,$j \in J$\, (giving \,$\sum
\limits_{j \in J} \beta_{j} = 1$),\, summing, and rearranging, we
obtain \,$A\bar{x} \leq b - A\sum \limits_{j \in J}
\beta_{j}p_{j}$.\, \smallskip In other words, it holds
\,$A(\bar{x} + \sum \limits_{j \in J} \beta_{j}p_{j}) \leq b$\,
subject to \,$\beta_{j}\geq 0,$\, \,$j \in J$,\, $\sum \limits_{j
\in J}\beta_{j} = 1$.\, By definition, the set \,$M$\, as the
convex hull of the vectors \,$p_{j}, j \in J$\,
 consists of all their convex combinations, so it
is convex. More precisely, all the points \,$y \in \Psi$\, are
specified by some possible convex combinations \,$\sum \limits_{j
\in J} \beta_{j}p_{j}$.\,  This allows us to conclude that the
following implication holds
$$A(\bar{x} + y) \leq b,\, \,\forall y \in \Psi\, \,\Rightarrow
\bar{x} \in \Phi - M.$$ Taking into account the arbitrary
selection of \,$\bar{x} \in \Psi_{1},$\, we thereby get the
desired inclusion \,$\Psi_{1} \subseteq  \Phi - M$.\,

Part II.\, Now, there is no loss of generality in taking
 some point \,$\bar{x}$\, ar\-bit\-ra\-ri\-ly from \,$\Phi - M$.\,
We aim at showing that \,$\bar{x} \in \Psi_{1}$.\, In this case,
the proof is trivial. Indeed, by our assumption, there is
fulfilled \,$A(\bar{x} + y) \leq b,\, \,\forall y \in \Psi$.\,
Since \,$p_{j} \in \Psi,$\, \,$\forall j \in J$,\, we obviously
have \,$A(\bar{x} + p_{j}) \leq b$\,  for all \,$j \in J$\, and
given points \,$p_{j} \in \mathbb{R}^{n}$,\, i.e. \,$\bar{x} \in
\Psi_{1}$.\, Due to our arbitrary choice of \,$\bar{x} \in \Phi -
M$,\, this therefore implies that \,$\Phi - M \subseteq
\Psi_{1}.$\, $\square$

 The right-hand side for a system of inequalities in (\ref{Gabid_Psi1}),
taken with all possible choices \,$j \in J$,\, characterized this
system as overdetermined. In what follows,  having reduced the
number of inequality constraints in (\ref{Gabid_Psi1}), we will
come latter to the refined representation of the Minkowski
difference \,$\Phi - M$.\,

For all \,$j \in J$,\, let us use the following notation:
\,$Ap_{j} = \tilde{b}_{j}$.\, Then we have
$$b^{i} - \hat{b}^{i} = b^{i} - \max \limits_{j \in J} \tilde{b}^{i}_{j}
\leq b^{i} - \tilde{b}^{i}_{j},\, \forall i \in K,\, \forall j \in
J.$$  For the system from (\ref{Gabid_Psi1}) can then be
formulated some subsystem of the form \,$Ax \leq s,$\, where

\centerline{\,$s \in \mathbb{R}^{r},$\,  $s^{T} = b^{T} -
\hat{b}^{T} = b^{T} - (\max \limits_{j \in J} \tilde{b}^{1}_{j},\,
\max \limits_{j \in J}\tilde{b}^{2}_{j},\, \ldots,\, \max
\limits_{j \in J} \tilde{b}^{r}_{j})$.\,} Define the set with the
following constraint representation
$$\Psi_{2} = \{ x \in \mathbb{R}^{n}: Ax \leq s\}.$$
\begin{theorem}\label{Gabid_Psi2-phi M} (Refined Minkowski difference for convex polyhedra
having different representations) Let be given an arbitrary
collection of vectors \,$p_{j} \in
   \mathbb{R}^{n}$,\, $j \in J,$
 \,$\Phi$\, be defined by (\ref{Gabid_Phi_lin_constr}),\, \,$\Phi
\neq \emptyset$,\, \,$M = conv\{p_{j}\}_{j \in J}$,\, then
$$\Phi - M = \Psi_{2}.$$
\end{theorem}
\proof  The core of the proof consists in showing that \,$\Psi_{1}
= \Psi_{2}$.\\
 Part I.\, There is no loss of generality in
selecting some fixed point \,$\bar{x}$\, from \,$\Psi_{2}$.\,
Writing
\\\centerline{$A\bar{x} \leq s = b - \hat{b} \leq b - \tilde{b}_{j} = b
- Ap_{j},\, \forall j \in J$,\,}\\ we see that \,$\bar{x}$\, lies
in \,$\Psi_{1}$.\, The arbitrariness of choosing \,$\bar{x} \in
\Psi_{2}$\, confirms that there is true the following inclusion:
\,$\Psi_{2} \subseteq \Psi_{1}$.\\
Part II.\, To prove the converse inclusion, take some point
\,$\bar{x}$\, arbitrarily from \,$\Psi_{1}$.\, Due to the special
construction of the right-hand side of constraint system from the
description of \,$\Psi_{1}$,\, for each \,$i \in K$,\, it is not
hard to see that the following subsystem of linear constraints
\begin{eqnarray}\nonumber
\sum \limits_{k = 1}^{n} a_{ik}\bar{x}^{k} \leq b^{i} - \sum
\limits_{k = 1}^{n} a_{ik}p_{1}^{k}, \\ \ldots \,\, \ldots \,\,
\ldots \,\, \ldots \ldots \,\, \ldots \,\, \ldots \nonumber \\
 \sum \limits_{k = 1}^{n} a_{ik}\bar{x}^{k} \leq b^{i} - \sum
\limits_{k = 1}^{n} a_{ik}p_{m}^{k} \nonumber
\end{eqnarray}
can be replaced by a unique inequality
$$\sum \limits_{k = 1}^{n} a_{ik}\bar{x}^{k} \leq b^{i} - \max \limits_{j \in J}\sum
\limits_{k = 1}^{n} a_{ik}p_{j}^{k}.$$ The truth of the previous
assertion is quite obvious since it holds \smallskip \\
\centerline{$\min \limits_{j \in J} (b^{i} - \sum \limits_{k =
1}^{n} a_{ik}p_{j}^{k}) = b^{i} - \max \limits_{j \in J}\sum
\limits_{k = 1}^{n} a_{ik}p_{j}^{k}$.}\smallskip \\ In the above
chain of the formulas, \,$\bar{x}^{k}$,\, $p_{j}^{k}$\, are the
$k-$th components of \,$\bar{x}$\, and \,$p_{j}$,\, respectively.
Besides, \,$a_{ik}$\, denotes  the $k-$th  element in the row
indexed by \,$i$\ of the matrix \,$A$,\, \,$b^{i}$\, corresponds
to the \,$i-$th item of the vector \,$b$.\,
 In consequence, the system
\,$A\bar{x}  \leq b - Ap_{j}$,\, $j \in J$\, can be reduced to the
 system \,$A\bar{x}  \leq s $\, with the much fewer number of constraints
but with the same dimension, i.e. \,$\bar{x} \in \Psi_{2}$.\, Due
to an arbitrary manner of selecting \,$\bar{x}$\, from
\,$\Psi_{1}$,\, this justifies the inclusion \,$\Psi_{1} \subseteq
\Psi_{2}$.\,  Through the forward and backward inclusions,
 we have \,$\Psi_{1} = \Psi_{2}$.\,  Consequently, from Theorem \ref{Gabid_Phi-M},
 it follows that \,$\Phi - M = \Psi_{2}$.\, $\square$

 In light of these assertions, inter alia, we
can now apply the formula of the Minkowski difference \,$\Phi -
 M$\, in solving the separation problem for \,$\Phi$\, and
 $M$.\,
                        To linearly separate the sets \,$\Phi$\,
                        and \,$M$,\, for a beginning, we thereby need to solve the
                        prob\-lem~(\ref{G_maxtfc}). Of course,
                        in the setting of (\ref{G_maxtfc}), we
                        take the set \,$\Phi - M$\, instead  of \,$\tilde{\Phi}$.\,
                        After that, we are now in the position of being able to analyze
                        the values of \,$t_{\Phi - M}(c^{*})$\, and \,$c^{*}$.\, This
                        analysis allows us to charac\-te\-ri\-ze  \,$\Phi$\, and \,$M$\, as
                        linearly (strongly or not) separable or
                        inseparable. In the separable case,
                        \,$c^{*}$\, represents the normal vector of
                        the best linear separator  for \,$\Phi$\, and
                        \,$M$\, with the maximal thickness. In
                        the case of the sets inseparability,
                        \,$c^{*}$\, corresponds to the best
                        pseudo-separator with the minimal
                        thickness. In other words, for
                        any case, the special setting of (\ref{G_maxtfc})
                        provides an opportunity for obtaining the optimal
                        thickness of the margin between
                        the supporting hyperplanes to the sets.
                        For the details relatively the terms of
                        separator, pseudo-separator, their thickness, and not only, the
                        interested reader is again directed to
                        \cite{G_Gabid6} (see \,p.~161).
Notice that if we need only inspect the issue of whether the
Minkowski difference \,$\Phi - M$\, fails to be strongly separable
from the origin of \,$\mathbb{R}^{n}$,\, then we can simply test
whether or not the origin satisfies the system of constraints
describing \,$\Phi - M$.\, Another question may now arise. How we
can reveal that the sets \,$\Phi$\, and  \,$M$\, are non-strongly
linearly separable from each other? Due to constraint structure of
\,$\Phi - M$,\, it is also not hard to check the question of
whether the origin of \,$\mathbb{R}^{n}$\, is the boundary point
of \,$\Phi - M$.\, Indeed, if some constraints is fulfilled at the
origin as equations, whereas all the others hold as the strict
inequalities, obviously this means then that the origin belongs to
the boundary of \,$\Phi - M$.\,

                        Let us note that the problem of computing
                        the distance between \,$\Phi$\, and
                        \,$M$\, can be solved by means of
                        projecting the origin of \,$\mathbb{R}^{n}$\, onto \,$\Phi - M$.\,
                        It will especially be expedient to apply
                        such a reduction of the problems in the
                        case when all of the right-hand side components \,$s_{k}$\, for the
                        system identifying the refined version of
                        the Minkowski difference \,$\Phi - M$\,
                        are negative. Since in this event, as it
                        was proved in \cite{Gabid_gab5},\, a projection problem
                        for convex polyhedron given by a system of linear
                        constraints can be solved by reduction to
                        the problem of projecting the origin onto
                        the polyhedron described by a finite collection of
                        points from \,$\mathbb{R}^{n}$.\, In its
                        own turn, this projection problem may be
                        solved utilizing one of the various problem settings
                        presented in \cite{Gabid_gab4}.\, The
                        proposed reduction makes wider a range of
                        suitable optimization tools which can effectively
                        be operated for solving the projection
                        problem in consideration.
 \subsection{Both Operands with a Half-space Representation}\label{Gabid_subsec2.4}\,
The purpose of this subsection is to obtain the description of the
Minkowski difference for the case where both operands have the
same setting in the form of the half-spaces intersection (or,
briefly, have a so-called half-space representation).

    Let us be given the two following polyhedra described as
    the intersection of closed half-spaces of \,$\mathbb{R}^{n}$\,:
    $$\Phi = \{ x \in \mathbb{R}^{n}: A_{1}x \leq b_{1}\},$$
$$\Psi = \{ x \in \mathbb{R}^{n}: A_{2}x \leq b_{2}\},$$
where \,$A_{1} \in \mathbb{R}^{r_{1} \times n}$,\, \,$A_{2} \in
\mathbb{R}^{r_{2} \times n}$\, are some nonvacuous matrices,\,
\,$b_{1} \in \mathbb{R}^{r_{1}}$,\, \,$b_{2} \in
\mathbb{R}^{r_{2}}$.\, By the same arguments already applied in
Sub\-sec\-ti\-on~2.3,\, there can be shown that \,$\Phi$\, and
\,$\Psi$\, are closed and convex sets. According to
Theo\-rem~\ref{G_Minkowski_1}, the Minkowski difference of sets
\,$\Phi$\,
 and \,$\Psi$\, can be expressed as follows
 $$\Phi - \Psi = \{ x \in \mathbb{R}^{n}: A_{1}(x+y) \leq b_{1},\, y \in \Psi\}.$$
 If we replace further the abstract constraint by the system defining \,$\Psi$,\,
  then we immediately obtain the following system of linear constraints:
  \begin{eqnarray}\nonumber
A_{1}x + A_{1}y \,\leq b_{1}, \\
A_{2}y \leq b_{2}.\nonumber
\end{eqnarray}
 This system can be equivalently rewritten in matrix
 form as follows:
  $$\Phi - \Psi = \{ z \in \mathbb{R}^{2n}: Dz \leq b\},$$
  where \,$D \in \mathbb{R}^{(r_{1} + r_{2}) \times 2n}$\, is a block-structured
  matrix,

 \centerline{$D = \left(
\begin{array}{cc}
   \enspace \,A_{1} & \mid A_{1} \\
   \hline
   \enspace  \Theta & \mid   A_{2}  \smallskip
  \end{array}
\right).$}\, Besides, the right-hand side of the system and the
vector of va\-ri\-ables have also the block structure, i.e. \,$b
\in \mathbb{R}^{r_{1} + r_{2}}$\, and
 \,$b = \left(
\begin{array}{c}
   \enspace \,b_{1}  \\
   \hline
   \enspace \,b_{2}
  \end{array}
\right),$\,  \,$z \in\mathbb{R}^{2n}$\,  and \,$z = (x \,\mid\,
y)$.\, Here, \,$\Theta \in \mathbb{R}^{r_{2} \times n}$\, denotes
the null matrix. Being the intersection of the closed half-spaces
in \,$\mathbb{R}^{2n}$,\, the Minkowski difference \,$\Phi -
\Psi$\, is the closed and convex point set.

 For computing the Euclidean distance between the sets  \,$\Phi$\, and \,$\Psi$,\,
  we can formulate and solve the following problem of minimizing the strongly convex
  quadratic function subject to the linear inequalities:
 \begin{eqnarray}
  \min \|z\|^{2} \label{Gabid_z^2} \\
  Dz \,\leq b,  \label{Gabid_Dz_b}
\end{eqnarray}
where the number of constraints is equal to \,$r_{1} + r_{2}$,\,
the amount of variables equals \,$2n$.\, For solving the same
problem of measuring the distance between \,$\Phi$\, and
\,$\Psi$,\, in \cite{G_Eremin},\, they  dealt with the following
optimization problem:
 \begin{eqnarray}
 \min \|x - y\|^{2} \label{Gabid_Er1}\\
A_{1}x \leq b_{1},  \label{Gabid_Er2}\\
A_{2}y \leq b_{2}. \label{Gabid_Er3}
\end{eqnarray}
Observe that this problem has the exactly similar number of
constraints and variables as the problem
(\ref{Gabid_z^2})--(\ref{Gabid_Dz_b}). Nevertheless, the objective
function of (\ref{Gabid_Er1})--(\ref{Gabid_Er3}) does not have the
strong convexity property. For this reason, the program
(\ref{Gabid_z^2})--(\ref{Gabid_Dz_b}) compares favorably with
(\ref{Gabid_Er1})--(\ref{Gabid_Er3}).

\subsection{Applications in Variational Inequalities Problems Related to a Concept of Linear Separability}
\label{Gabid_subsec2.5}\,
 Naturally, the analytical representation
of the Minkowski difference of sets already has its own utility
and independent significant applications in various fields of
mathema\-ti\-cal sciences. Likewise, the Minkowski difference
operation can take its place now as a~mathematical tool ready for
new applications. In regard the topic of interest, it is quite
clear that the Minkowski difference is a~tool suitable for dealing
with solving the variational inequalities that are closely
relevant to a concept of the linear separability of sets.

 Let us consider the following variational inequalities, which consist in determining
a~non\-zero vector \,$c \in \mathbb{R}^{n}$\, such that
            \begin{equation}\label{Gab_VAR_1}
\langle c, x - y \rangle \geq \Delta, \quad x \in A,\, y \in B,\,
\Delta > 0,
            \end{equation}
            \begin{equation}\label{Gab VAR_2}
\langle c, x - y -c \rangle \geq 0, \quad x \in A,\, y \in B,
              \end{equation}
              \begin{equation}\label{Gab VAR_3}
\langle c, x - y \rangle \geq 0, \quad x \in A,\, y \in B.
              \end{equation}
              Note that the first two of these variational inequalities correspond
              to the strong linear separability term in a sense that the inequalities
              are solvable in the case when the sets $A$\, and \,$B$\,
              are strongly linearly separable.
               The third inequality is
              closely connected with the potentially non-strong linear
              separability of the considered sets.
              Clearly, using the Minkowski difference \,$A - B$,\,  for the above-mentioned inequalities
              (\ref{Gab_VAR_1})--(\ref{Gab VAR_3}), a~re\-duc\-ti\-on
              can be made to the following
              variational inequalities, respectively:
             \begin{equation}\label{Gab VAR_4}
\langle c, z \rangle \geq \Delta, \quad x \in A - B,\, \Delta > 0,
              \end{equation}
               \begin{equation}\label{Gab VAR_5}
\langle c, z - c \rangle \geq 0, \quad z \in A - B,
              \end{equation}
               \begin{equation}\label{Gab VAR_6}
\langle c, z \rangle \geq 0, \quad z \in A - B.
              \end{equation}
                             Of course, the goal is to find the nonzero vector \,$c \in
              \mathbb{R}^{n}$\, satisfying (\ref{Gab VAR_4})--(\ref{Gab VAR_6}).\,
              A set of possible solutions for (\ref{Gab VAR_6}) coincides
              with \,$W_{A - B}\ \backslash \{{\bf 0}\}$.\, For the more general than
              convex polyhedral setting, the proof of this fact
              was represented in \cite{G_Gabid7}(see p. 712).
              Moreover, it was proved that the inequality (\ref{Gab VAR_6})
              has nonzero solutions if and only if \,$W_{A - B}\ \neq \{{\bf
              0}\}$\, or, equivalently, \,${\bf 0} \notin int(A -
              B)$.\, For the convex sets \,$A$\,  and \,$B$\, with nonempty
              interiors, due to Lemma 3.16\, from \cite{G_Gabid7}, it holds
$$int(A - B) = int(A) - int(B).$$ This means that (\ref{Gab VAR_6})\,
(in tandem with (\ref{Gab VAR_3})) is solvable if and only if the
convex polyhedra \,$A$\, and \,$B$\, with nonempty interiors have
no the common interior points. For convex and closed sets \,$A$\,
and \,$B$,\, a solution set of (\ref{Gab VAR_4})\, coincides with
\,$V_{A - B}$.\, In \cite{G_Gabid7},\, it was justified that the
variational inequality (\ref{Gab VAR_4})\, together with
(\ref{Gab_VAR_1})\, has solutions if and only if \,$V_{A - B} \neq
\emptyset$\, or, equivalently, \,${\bf 0} \notin A - B$.\, Thus,
the absence of any common points for \,$A$\, and \,$B$\, is the
 condition for the solvability of (\ref{Gab VAR_4})\, and
(\ref{Gab_VAR_1}).\, Formally, the fulfillment of this condition
can be verified with the help of projecting the origin of
\,$\mathbb{R}^{n}$\, onto \,$A - B$.\, If as the result of
projecting we obtain \,${\bf P}_{A - B}({\bf 0}) \neq {\bf 0}$,\,
then \,$c = {\bf P}_{A - B}({\bf 0})$\, represents the solution
for (\ref{Gab VAR_4}).\, Otherwise, the inequalities (\ref{Gab
VAR_4})\, and (\ref{Gab_VAR_1})\, have no solutions.
               For the variational inequality which consists in
              determining a~vec\-tor \,$c \in \mathbb{R}^{n}\backslash \{{\bf
              0}\}$\, satisfying (\ref{Gab VAR_5}), a solution set
              coincides with \,$\Omega_{A - B}\backslash \{{\bf
              0}\}$.\, Let us remind that the considered pair of convex
               polyhedra \,$A$\, and \,$B$\, are closed.  Then the boundedness of one of these
               sets implies the closedness of \,$A - B$.\,
               By the closedness and convexity properties of \,$A - B$,\,
               due to Theorem 3.3\, from \cite{G_Gabid6},\, we
                immediately obtain \,${\bf P}_{A - B}({\bf 0}) \in Bd(\Omega_{A - B})$.\,
               Having obtained \,${\bf P}_{A - B}({\bf 0}) = {\bf
               0}$,\, we can conclude that \,$\Omega_{A - B} = \{{\bf 0}\}$,\,
                since Theorem 3.2\, from \cite{G_Gabid6}\, yields
                \,$$\min \limits_{x \in A- B} \|x\|^{2} =
                \max \limits_{y \in \Omega_{A- B}} \|y\|^{2}.$$\,
                The fulfillment of \,$\Omega_{A - B} = \{{\bf
                0}\}$\, is equivalent to having \,${\bf 0} \in A -
                B$,\, since \,$E(\Omega_{A - B})\backslash\{{\bf 0}\} = V_{A - B}.$\,

              For the Minkowski difference operation, we underline the next
application which consists in finding the nearest points of sets
$A$\, and \,$B$\, by solving the system
 \begin{equation*}
\left\{\begin{aligned}
\langle c, x - \bar{x} \rangle \geq 0, \quad x \in A,\\
\langle c, \bar{y} - y \rangle \geq 0, \quad y \in B,
\end{aligned}\right.
\end{equation*}
where \,$c = {\bf P}_{A - B}({\bf 0})$\, is the projection of the
origin of \,$\mathbb{R}^{n}$\, onto \,$A - B$.\, Let us note that
the projection problem can be reduced to the maximin problem of
the type (\ref{G_maxtfc}). To solve the maximin problem, we can
  apply some software package after a~simple transition to the minimax
  problem as follows
  \begin{multline}\nonumber
\max \limits_{w \in \chi} \,\,\inf \limits_{x \in A - B} \,\langle
x,w \rangle = -\min \limits_{w \in \chi} \,\,(-\inf \limits_{x \in
A - B} \,\langle x,w \rangle)\, \,= \\ = -\min \limits_{w \in
-\chi} \,\,\sup \limits_{x \in A - B} \,\langle x,w \rangle =
-\min \limits_{w \in \chi} \,\,\sup \limits_{x \in A - B}
\,\langle x,w \rangle, \text{where}\\  \,\,\chi := \{c \in
\mathbb{R}^{n}: \|c\| = 1\},\,\,-\chi:= \{-x,\, x \in \chi\}.\,
\end{multline}
We notice that, for instance, a package Optimization Toolbox in
MATLAB contains the function Fminimax which is usable for solving
the minimax constraint problem.
\section{Conclusions}
\label{Gabid_subsec2.6}\, We have presented the novel analytical
representation of the Minkowski difference for convex polyhedra
given by different ways. In particular, we have considered the
following cases where:
\begin{itemize}
    \item both operands under the Minkowski difference operation
are similarly determined as convex hulls of finite collections of
points,
 \item the operands have the different performance (or more
precisely,  the first operand is expressible by a linear
constraint system, and the second one is given in terms of convex
hull), \item  both operands have the same representation as the
intersection of the closed half-spaces.
    \end{itemize}
   It can be concluded that thanks to the obtained results, there
appeared the possibility to investigate the relevant problems such
as: problem of linear se\-pa\-ra\-tion of the convex polyhedra in
the Euclidean space, variational inequalities problems, problem of
finding the distance between the convex polyhedra by projecting
the origin of the Euclidean space onto a convex polyhedron, the
problem of determining the closest points of convex polyhedra.

\end{document}